\documentclass{article}
\usepackage{theorem}
\theorembodyfont{\upshape}
\usepackage{latexsym}
\usepackage{amssymb,amsmath,dsfont}
\usepackage{graphicx,float}
\usepackage{subfigure}
\usepackage{url}

\newtheorem{thm}{Theorem}[section]

\def\dsp{\displaystyle}
\def\abs#1{\vert #1 \vert}
\def\jump#1{\left[ #1 \right]}
\newcommand{\R}{\mathbb{R}}
\newcommand{\Z}{\mathbb{Z}}
\newcommand{\WW}{\textbf{W}}
\newcommand{\BB}{\textbf{B}}
\newcommand{\U}{\textbf{U}}
\newcommand{\F}{\textbf{F}}
\newcommand{\fig}{F{\sc{ig}}.~}
\newcommand{\SE}{\mathcal{\textbf{S}}}
\def\vecdeux#1#2{\left(\begin{array}{cc}
           #1\\
           #2
          \end{array}
          \right)
         }
\begin{document}
\title{A kinetic scheme for transient mixed flows in non uniform closed
pipes: a global manner to upwind all the source terms}

\author {C. Bourdarias, M. Ersoy, S. Gerbi \\
Universit\'{e} de Savoie, Laboratoire de Math\'ematiques, \\
73376 Le Bourget-du-Lac Cedex,
France.}
\date{\today}

\maketitle

\begin{abstract} We present a numerical kinetic scheme for an unsteady mixed pressurized 
and free surface model. This model has a
source term depending on both the space variable and the unknown $U$ of the system.  
Using the  Finite Volume and Kinetic (FVK) framework, 
we propose  an approximation of the source terms following the principle of interfacial 
upwind with a kinetic interpretation.
Then, several numerical tests are presented.
\end{abstract}
\section{Introduction}
In this paper,  we study  a way to upwind the source terms of a mixed flows model in non uniform closed water pipes  
in a one
dimensional framework. 
In the case of free surface incompressible flows, the model is called \textbf{FS}-model and 
it is an extension of classical Saint-Venant models. 
When the pipe is full, we introduce the pressurized model, called \textbf{P}-model, 
which describes the evolution of
a compressible inviscid flow and is close to gas dynamics equations in a nozzle. In order to cope with the transition
 between a free surface and a pressurized model, we use a mixed model called \textbf{PFS}-model. It is based on balance laws 
and provides an hyperbolic system with source terms corresponding to the inclination of the pipe (seen as a topography term), 
the section variation, the curvature  and the friction.

Several ways to compute the numerical approximation of  conservation laws with source terms have already been
investigated by many authors. The main difficulty is to preserve numerically some properties
satisfied by the continuous model: 
the invariant domain, the well-balanced property for instance. The Finite Volume methods are largely used  since 
they present  the remarkable property to be domain invariant  (for instance, for Saint-Venant
equations, to be water height conservative). Some Well-Balanced Finite Volume  schemes,
 introduced by Greenberg \emph{et al} \cite{GL96}, preserve steady states.
All these methods are based on two principles:
firstly, the conservative quantities are cell-centered as usual finite volume schemes, and secondly the source 
terms are upwinded at the cell interfaces. 
 
In this paper, we consider a  particular Finite Volume-Kinetic scheme built to compute the numerical
solution of \textbf{PFS}-model. This scheme is  based on the
classical kinetic interpretation \cite{PS01} of the system. 

The source terms appearing in the \textbf{PFS}-model are either
conservative, non-conservative or else. All source terms are upwinded at the cell interfaces:
we use the definition of the DLM theory \cite{DLM95} to define the non-conservative products.
The particular case of the friction term which is neither conservative nor
non-conservative will be upwinded using the notion of \emph{dynamic slope}. 
The source terms are taken into account  in the numerical fluxes and are computed from the
microscopic ones, obtained through the concept of  potential barrier.

The paper is organized as follows. In  the second section, we describe the \textbf{PFS}-model \cite{BEG09_2} and 
focus on the  source terms. 
The detailed description of the method used to deal with the transition points (when a change of state occurs) 
is not presented (see \cite{BEG09_2} for more details on this topic). 
We state some theoretical
properties of the system. In the third section, we give the kinetic formulation of the \textbf{PFS}-model and
the corresponding kinetic scheme. 
In the fourth and last section, several numerical tests are provided.
{\small
\begin{itemize}
\subsubsection*{Notations concerning geometrical quantities}
\item $\theta(x)$:  angle of the inclination of the main pipe axis $z = Z(x)$ at position $x$
\item $\Omega(x) = S(x)$:  cross-section area of the pipe orthogonal to the axis $z=Z(x)$
\item $R(x)$:  radius of the cross-section $S(x)$ orthogonal to the axis $z=Z(x)$
\item $\Omega(t,x)$ : free surface cross-section area orthogonal to the axis $z=Z(x)$
\item $\sigma(x,z)$: width of the cross-section $\Omega$ at  position $x$ and altitude $z$
\subsubsection*{Notations concerning the \textbf{PFS}-model}
\item $p(t,x,y,z)$: pressure
\item $\rho_0$:  density of the water at atmospheric pressure $p_0$
\item $\rho(t,x,y,z)$:  density of the water at the current pressure
\item $\dsp\overline{\rho}(t,x) = \frac{1}{S(x)}\int_{\Omega(x)} \rho(t,x,y,z)\,dy\,dz$: 
mean value of $\rho$ over $\Omega(x)$
\item $c$: sonic speed
\item $\dsp A(t,x) =  \frac{\overline{\rho}(t,x)}{\rho_0} S(x)$:  equivalent wet area
\item $u(t,x)$: velocity
\item $Q(t,x) = A(t,x) u(t,x)$:  discharge
\item $E$: state indicator equal to $E=0$ if the flow is free surface, $E=1$ otherwise
\item $\SE$: the physical wet area equal to $A$ if the state is free surface, $S$ otherwise
\item $\mathcal{H}(\SE)$: the $Z$-coordinate of the water level equal to  $\mathcal{H}(\SE)=h(t,x)$ 
if the state is free surface, $R(x)$ otherwise
\item $p(x,A,E)$: mean pressure over $\Omega$
\item $K_s>0$: Strickler coefficient depending on the material 
\item $P_m(A)$: wet perimeter of $A$ (length of the part of the channel’s section in contact with the water) 
\item $R_h(A) = \dsp \frac{A}{P_m(A)}$: hydraulic radius
\item  Bold characters are used for vectors, except for $\SE$
\end{itemize}
}

\section{A model for  unsteady water flows in pipes }\label{SectionAModelForUnsteadyWaterFlowsPipes}
The \textbf{PFS}-model \cite{BEG09_2}  
is a mixed model of a pressurized (compressible) and free surface  (incompressible) flow in a one
dimensional  rigid pipe with variable cross-section. The pressurized parts of the flow correspond 
to  a full pipe whereas  the section is not completely filled for the free surface flow. 
The \textbf{F}ree \textbf{S}urface part of the model is derived  by writing the $3$D Euler
incompressible equations and by averaging over orthogonal sections to the privileged axis of the
flow. 
In the same spirit, by writing the Euler isentropic and compressible equations with the linearized
pressure law 
$\dsp p(t,x,y,z)=p_0 + \frac{1}{c^2}(\rho(t,x,y,z) - \rho_0)$, we obtain a Saint-Venant like
system of equations in the ``FS-equivalent'' variable 
$\dsp A(t,x) = \frac{\overline{\rho}(t,x)}{\rho_0} S(x)$, $Q(t,x)=A(t,x) u(t,x)$ 
which takes into account the compressible effects (for a detailed derivation, see \cite{BEG09_2}).

In order to deal with the transition points (that is, when a change of state occurs), we introduce a
state indicator variable $E$ which is equal to $1$
if the state is pressurized and to $0$ if the state is free surface. We define the \emph{physical wet area} by:
$$
\SE = \SE(A,E) = \left\{
\begin{array}{lll}
S & \textrm{ if } & E = 1 ,\\
A & \textrm{ if } & E = 0.
\end{array}
\right.
$$
The pressure law is given by a  mixed ``hydrostatic'' (for the free surface part of the flow) 
and ``acoustic''  type (for the pressurized part of the flow) as follows:
\begin{equation}\label{PFSPressureLaw}
 \dsp p(x,A,E) =  c^2(A-\SE) + gI_1(x,\SE)\cos\theta\,
\end{equation} where $g$ is the gravity constant, $c$ the sonic speed of the water (assumed to be
constant) and $\theta$ 
the inclination of the pipe. The term $I_1$ is the classical hydrostatic pressure: 
$$\dsp I_1(x,\SE) =
\int_{-R}^{\mathcal{H}(\SE)}(\mathcal{H}(\SE)-z) \sigma \,dz$$ where $\sigma(x,z)$ is the width of the
cross-section, $R=R(x)$ the radius of the cross-section and $\mathcal{H}(\SE)$ is the $z$-coordinate of the
free surface over the main  axis $Z(x)$.

\noindent The defined pressure (\ref{PFSPressureLaw}) is  continuous throughout the
transition points and we define the  \textbf{PFS}-model by:
\begin{equation}\label{PFS}
\left\{
\begin{array}{lll}
\partial_{t}(A) + \partial_{x}(Q) &=&0 \\
\partial_{t}(Q) + \partial_{x} \left(\dsp \frac{Q^2}{A}+p(x,A,E)\right) &=&\dsp-g\, A\,Z' + Pr(x,A,E) \\
 & &\dsp-G(x,A,E) \\
 & &\dsp - K(x,A,E) \dsp\frac{Q|Q|}{A}
\end{array}
\right.\,
\end{equation}
where $z=Z(x)$ is the altitude of the
main pipe axis. The terms $Pr$, $G$ and $K$ denote respectively
the pressure source term, a curvature term   and the friction:
$$\begin{array}{lll}
Pr(x,A,E) &=& \dsp c^2\left(\frac{A}{\SE}-1\right)\dsp S'+
g\,I_2(x,\SE)\cos\theta , \\
 G(x,A,E) &=& \dsp g\,A\,  \overline{Z}(x,\SE) = \dsp g\,A\,  \left(\mathcal H(\SE)-I_1(x,\SE)/\SE\right) \dsp
(\cos\theta)',
\\
K(x,A,E) &=& \dsp \frac{1}{K_s^{2} R_h(\SE)^{4/3}}\, 
  \end{array}
$$
where we have used the notation $f'$ to denote the  derivative with respect to the space variable $x$ 
of any function $f(x)$. 
The term $I_2$ is the hydrostatic pressure source term defined by:  $\dsp I_2(x,\SE) =
\int_{-R}^{\mathcal{H}(\SE)}(\mathcal{H}(\SE)-z) \partial_x\sigma \,dz\,.$ The term   
$K_s>0$ is the Strickler coefficient depending on the material and $R_h(\SE)$ is the hydraulic radius. 

The System (\ref{PFS}) has the following properties:
\begin{thm}\label{ThmPFSModel}
\begin{enumerate}
\item[]
\item System (\ref{PFS}) is strictly hyperbolic on $\left\{A(t,x)>0\right\}\,.$
\item For smooth solutions, the mean velocity $u = Q/A$ satisfies
\begin{equation}\label{ThmPFSEquationForU}
\begin{array}{c}
\partial_t u + \partial_x \left(\displaystyle\frac{u^2}{2} + c^2 \ln(A/{\SE}) + g\mathcal{H}(\SE)\cos\theta + gZ \right) \\ = -g K(x,A,E) u|u|
\leqslant 0.
\end{array}
\end{equation}
\item The still water steady state, for $u = 0$, reads:
\begin{equation}\label{ThmPFSSteadyState}
c^2 \ln(A/{\SE}) + g\mathcal{H}(\SE)\cos\theta + gZ = 0.
\end{equation}
\item System (\ref{PFS}) admits a mathematical entropy $$\mathcal{E}(A,Q,E) =\displaystyle \frac{Q^2}{2A} + c^2 A \ln(A/{\SE})+ c^2 S + g A
\overline{Z}(x,\SE)\cos\theta + gAZ $$
\noindent which satisfies the entropy relation for smooth solutions
\begin{equation}\label{ThmPFSEntropy}
\partial_t \mathcal{E} +\partial_x \big((\mathcal{E}+p(x,A,E))U\big) = -gAK(x,A,E) u^2 |u| \leqslant 0\,.
\end{equation}
\end{enumerate}
\end{thm}
In what follows, when no confusion is possible, the term $K(x,A,E)$ will be noted simply $K(x,A)$
for free surface states  and $K(x,S)$ for pressurized states.
\section{The Kinetic approach}\label{SectionKineticInterpretationPFSModel}
The kinetic formulation (\ref{ThmKineticFormulationPFS}) is a (non physical) microscopic
description of the \textbf{PFS}-model
provided by a given real function $\chi:\R\to\R$ satisfying the following  properties:
\begin{equation*}
\chi(\omega)=\chi(-\omega) \geqslant 0\;,\;
\int_{\R} \chi(\omega) d\omega =1,
\int_{\R} \omega^2 \chi(\omega) d\omega=1 .
\end{equation*} 
It permits to  define  the density of particles, by a so-called \emph{Gibbs equilibrium}, 
$\dsp
\mathcal{M}(t,x,\xi) =
\frac{A(t,x)}{b(t,x)} \chi\left(\frac{\xi-u(t,x)}{b(t,x)}\right)$ where $b(t,x) = b(x,A(t,x),E(t,x))$ with 
\begin{equation*}
b(x,A,E)= \left\{
\begin{array}{lll}
\dsp\sqrt{g\,\frac{I_1(x,A)}{A}\cos\theta} & \textrm{ if } & E = 0,\\
\dsp\sqrt{g\,\frac{I_1(x,S)}{A}\cos\theta+c^2} & \textrm{ if } & E = 1.\\
\end{array}
\right.
\end{equation*}
\subsection{The mathematical kinetic formulation}
The Gibbs equilibrium $\mathcal{M}$ is related 
to the \textbf{PFS}-model by the classical 
 \emph{macro-micro}scopic kinetic relations:
\begin{equation}\label{PFSWetAreaByMacroMicroRelations}
A  = \dsp\int_{\R} \mathcal{M}(t,x,\xi)\,d\xi\,,
\end{equation}
\begin{equation}\label{PFSDischargeByMacroMicroRelations}
Q  = \dsp\int_{\R} \xi\mathcal{M}(t,x,\xi)\,d\xi\,,
\end{equation}
\begin{equation}\label{PFSFluxByMacroMicroRelations}
\dsp\frac{Q ^2}{A }+A\,b(x,A,E)^2  = \dsp\int_{\R} \xi^2 \mathcal{M}(t,x,\xi)\,d\xi\,.
\end{equation}
From the relations \eqref{PFSWetAreaByMacroMicroRelations}--\eqref{PFSFluxByMacroMicroRelations}, 
the non-linear \textbf{PFS}-model can be viewed as a single linear equation involving the non-linear quantity 
$\mathcal{M}$:
\begin{thm}[Kinetic Formulation of the \textbf{PFS}-model]\label{ThmKineticFormulationPFS}
$(A,Q)$ is a strong solution of
System (\ref{PFS}) if and only if ${\mathcal{M}}$ satisfies the kinetic transport equation:
\begin{equation}\label{KineticFormulationPFS}
\partial_t \mathcal{M}+\xi \cdot \partial_x\mathcal{M} - g\phi 
\,\partial_\xi \mathcal{M} = \mathcal{K}(t,x,\xi)
\end{equation}
for a collision term $\mathcal{K}(t,x,\xi)$ which satisfies for $(t,x)$ a.e. 
$$\dsp  \int_{\R} \vecdeux{1}{\xi} \mathcal{K}(t,x,\xi)\,d\xi = 0.$$
\end{thm}
The source terms are defined as:
\begin{equation}\label{PFSSourceTermPhi}
\phi(x,\WW) = \BB(x,\WW)\cdot\partial_x \WW
\end{equation}
with \begin{equation}\label{WW}
\dsp\WW = \left(Z+\int_x K(x,A) u\abs{u}\,dx,\;S,\;\cos\theta\right)      
     \end{equation}
\noindent and
$
\BB =
\left\{
\begin{array}{ll}
\dsp \left(1,\;-\frac{c^2}{g}\left(\frac{A-S}{A\,S}\right)-\frac{\gamma(x,S)\cos\theta}{A},\;\overline{Z}(x,S)\right) & \textrm{ if  } E = 1,\\
\dsp \left(1,\;-\frac{\gamma(x,A)\cos\theta}{A},\;\overline{Z}(x,A)\right) & \textrm{ if  } E = 0\\
\end{array}
\right.
$

\noindent where $ I_2(x,\SE) $ reads $\gamma(x,\SE)
S'$ for some function $\gamma$ (depending on the geometry of the pipe).

We call the term $\dsp \frac{d}{dx}\left( Z+\int_x K(x,A) u\abs{u}\,dx\right)$ the \emph{dynamic slope} 
since it is time and space variable dependent. 
\subsection{The kinetic scheme}
Based on the above kinetic formulation \eqref{KineticFormulationPFS}, 
we construct easily a Finite Volume scheme where 
the source terms are upwinded by a generalized kinetic scheme with reflections \cite{PS01}. 

To this end, let us
consider a uniform mesh on $\R$ where  cells are denoted for every $i\in \Z$ by $m_i = (x_{i-1/2},x_{i+1/2}),$ with  $x_i
=\dsp\frac{x_{i-1/2}+x_{i+1/2}}{2}$ and $\Delta x =x_{i+1/2}-x_{i+1/2} $ the space-step. 
We  consider a time discretization $t^n$
defined by $t^{n+1}=t^n+\Delta t^n$ with $\Delta t^n$ the time-step.
We note $\U_i^n=(A_i^n,Q_i^n)$, $\dsp u_i^n = \frac{Q_i^n}{A_i^n}$, $\mathcal{M}_i^n$ 
the cell-centered approximation of $\U = (A,Q)$, $u$ and $\mathcal{M}$ on the cell $m_i$ at time $t^n$.

If $\dsp\WW$ is  $\dsp\left(Z+\int_x K(x,A) u\abs{u}\,dx,\;S,\;\cos\theta\right)$,  
its piecewise constant representation is given by,   
$\WW(t,x) = \dsp \WW_i(t) \mathds{1}_{m_i}(x)$ 
where $\WW_i(t)$ is defined as $\WW_i(t) = \dsp \frac{1}{\Delta x}\int_{m_i}\WW(t,x)\,dx$ for instance.

\noindent Denoting by $\WW_i$ and $\WW_{i+1}$ the left and the right states of the cell interface $x_{i+1/2}$, 
and using the ``straight lines'' paths (see  \cite{DLM95}) 
$$\Psi(s,\WW_i,\WW_{i+1}) = s\WW_{i+1}+(1-s)\WW_i,\,s\in[0,1],$$ 
we define the non-conservative product $\phi(t,x_{i+1/2})$ by writing:
\begin{equation}\label{DefinitionOfNonConservativeProductPhi}
\jump{\WW}(t)\cdot\int_0^1 \BB\left(t,\Psi(s,\WW_i(t),\WW_{i+1}(t))\right)ds
\end{equation}
where  $\jump{\WW}(t):=\WW_{i+1}(t)-\WW_i(t)$, is the jump of $\WW(t)$ across the discontinuity localized 
at $x=x_{i+1/2}$. 
As the first component of $\BB$
is $1$, we recover  the classical interfacial upwinding for the  term $Z$ 
(appearing e.g. in  Saint-Venant equations) since it is a
conservative term.

Neglecting the collision kernel \cite{PS01} and using the fact that 
$\phi = 0$ on the cell $m_i$ (since $\jump{\WW} \equiv 0$), the kinetic transport equation
(\ref{KineticFormulationPFS}) simply reads:
\begin{equation}
\left\{\begin{array}{l}
\dsp
\frac{\partial}{\partial t} f+\xi \cdot \frac{\partial}{\partial x}
f = 0\\ \\
f(t_n,x,\xi) = \mathcal{M}(t_n,x,\xi)
\end{array}\right.,(t,x,\xi)\in [t_n,t_{n+1})\times m_i\times\R
\end{equation}
and thus it may be discretized as follows:
\begin{equation}\label{DiscretisationOfKineticTransportEquation}
f_i^{n+1}(\xi) = \mathcal{M}_i^n(\xi)+\frac{\Delta t^n}{\Delta x}\,
\xi \, (\mathcal{M}_{i+\frac{1}{2}}^-(\xi)-\mathcal{M}_{i-\frac{1}{2}}^+(\xi))
\end{equation}
where the contribution of the source term  
is included into the microscopic numerical  fluxes $\mathcal{M}^{\pm}_{i\pm 1/2}$. 
This is the
principle of interfacial  source upwind. Using the \emph{macro-micro}scopic relations
\eqref{PFSWetAreaByMacroMicroRelations}--\eqref{PFSFluxByMacroMicroRelations} and integrating Equation
(\ref{DiscretisationOfKineticTransportEquation}) against $\xi$ and $\xi^2$, we obtain the Finite Volume scheme:
\begin{equation}\label{FVScheme}
\U_i^{n+1} = \U_i^n+\frac{\Delta t^n}{\Delta
x}\,(\F_{i+\frac{1}{2}}^{-}- \F_{i-\frac{1}{2}}^+)
\end{equation}
where the numerical fluxes are computed by :
\begin{equation}
\dsp \F_{i+\frac{1}{2}}^\pm = 
\int_{\R} \left( \begin{array}{c}
\xi \\
\xi^2
\end{array}\right)\,\mathcal{M}_{i+\frac{1}{2}}^\pm(\xi)\,d\xi\,.
\end{equation}
Following  \cite{PS01} (or \cite{BEG09_1}), the microscopic fluxes are given by:
\begin{equation}\label{MicroscopicInterfaceFluxes}
\begin{array}{lll}
\mathcal M_{i+1/2}^{-}(\xi) &=&  
\overbrace{\dsp \mathds{1}_{\{\xi>0\}}\mathcal M_i^n(\xi)}^{{\textrm{positive  transmission}}}
+\overbrace{\mathds{1}_{\{\xi<0,\xi^2-2g\Delta\phi^n_{i+1/2}<0\}}\mathcal M_i^n(-\xi)}^{\textrm{reflection}}\\ 
&+& \underbrace{\dsp \mathds{1}_{\{\xi<0,\xi^2-2g\Delta\phi^n_{i+1/2}>0\}}
\mathcal M_{i+1}^n\left(\dsp-\sqrt{\xi^2-2g\Delta  \phi^n_{i+1/2}}\right)}_{\textrm{negative   transmission}},\\
 & & \\
\mathcal M_{i+1/2}^{+}(\xi) &=& \overbrace{\dsp \mathds{1}_{\{\xi<0\}}\mathcal M_{i+1}^n(\xi)}^{\textrm{negative transmission}}
+\overbrace{\mathds{1}_{\{\xi>0,\xi^2+2g\Delta\phi^n_{i+1/2}<0\}}\mathcal M_{i+1}^n(-\xi)}^{\textrm{reflection}}\\ 
&+& \underbrace{ \dsp \mathds{1}_{\{\xi>0,\xi^2+2g\Delta\phi^n_{i+1/2}>0\}}
\mathcal M_{i}^n\left(\dsp\sqrt{\xi^2+2g\Delta\phi^n_{i+1/2}}\right)}_{\textrm{positive  transmission}}.
\end{array}\,
\end{equation}
The term $\Delta\phi^n_{i\pm 1/2}$ in (\ref{MicroscopicInterfaceFluxes})  is the upwinded source term
(\ref{PFSSourceTermPhi}). It also plays the role of the potential bareer:  
the term $\xi^2\pm 2g\Delta \phi^n_{i+1/2}$ is the jump condition for a particle with a kinetic speed $\xi$  
which is necessary  to
\begin{itemize}
\item be reflected: this means that the particle has not enough kinetic energy $\xi^2/2$ to overpass the potential barrier (reflection in
(\ref{MicroscopicInterfaceFluxes}))),
\item overpass the potential barrier with a positive speed (positive  transmission in (\ref{MicroscopicInterfaceFluxes})),
\item overpass the potential barrier with a negative speed (negative  transmission in (\ref{MicroscopicInterfaceFluxes}))).
\end{itemize}
Taking an approximation of the non-conservative product $\phi$ 
(\ref{DefinitionOfNonConservativeProductPhi}), 
the potential barrier $\Delta  \phi^n_{i+ 1/2}$ has the following expression:
\begin{equation}\label{ApproximationSourceTerm}
\Delta\phi_{i+1/2}^n= \jump{\WW}(t_n)\cdot  \BB\left(t_n,\Psi\left(\frac{1}{2},\WW_i(t_n),\WW_{i+1}(t_n)\right)\right)
\end{equation}
\noindent Next, with the simplest choice of the $\chi$-function  $\dsp 
\chi(\omega)=\frac{1}{2\sqrt{3}}\mathds{1}_{[-\sqrt{3},\sqrt{3}]}(\omega),$
which allows to compute easily numerical fluxes, we have:
\begin{thm}\label{ThmStability}
\begin{enumerate}
\item[]
\item  Under the CFL condition $\dsp\frac{\Delta t^n}{\Delta x} \max_{i\in\Z}\left(|u_i^n| +\sqrt{3}c\right)< 1,$ 
the kinetic scheme \eqref{FVScheme}--\eqref{MicroscopicInterfaceFluxes} keeps $A$ positive, i.e. $A_i^n\geqslant 0$ 
if it is  initially true.
\item The kinetic scheme \eqref{FVScheme}--\eqref{MicroscopicInterfaceFluxes} allows to compute the drying and flooding area.
\end{enumerate}
\end{thm}
\section{Numerical results}
Let us recall that the zero water level corresponds to the main pipe axis. 
The piezometric head (or line) is defined  by:
$$\displaystyle  piezo = z + p \; \mbox{ with }
\left\{\begin{array}{l}
\displaystyle p = 2R+\frac{c^2 \, (A -S)}{g\,S}
\mbox{ if the flow is pressurized}\\
p = h \mbox{ if the flow is free surface},
\end{array}
\right.
$$
where $h$  is the water height.
\paragraph{Comparison with the VFRoe scheme \cite{BEG09_2}.\newline}
We compare the result obtained by the presented kinetic scheme with the upwinded VFRoe method \cite{BEG09_2}. 

The numerical experiment is performed in the case of an expanding  $5$ $m$ long closed
circular pipe  at altitude  $Z_0 = 1\,m$ with $0$ slope (slope of the main pipe axis). 
The upstream diameter is $2\;m$ and the downstream diameter is $2.2\;m$. 
The friction is not considered for the first test and is set to $0$. The
simulation starts with a still free surface steady state.
The upstream boundary condition  is a prescribed piezometric line (increasing linearly from 
$1\;m$ to $3.2\;m$ in $5\;s$) while the downstream discharge is kept constant 
to $0\,m^3/s$. 
The other parameters are $N=100$ (discretization points), CFL$=0.8$ and the sound speed is
$20\;m/s$.

\begin{figure}[H]
\begin{center}
\subfigure[ Discharge]
{
\includegraphics[scale = 0.78]{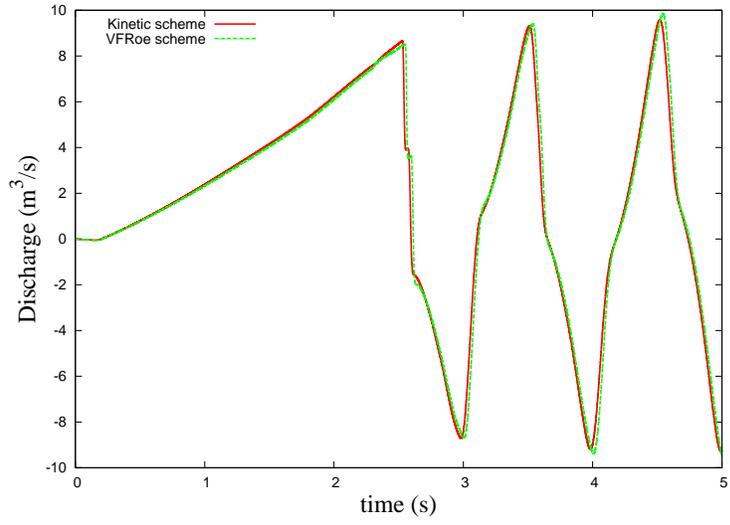}
}
\\
\subfigure[ Piezometric line]
{
  \includegraphics[scale = 0.78]{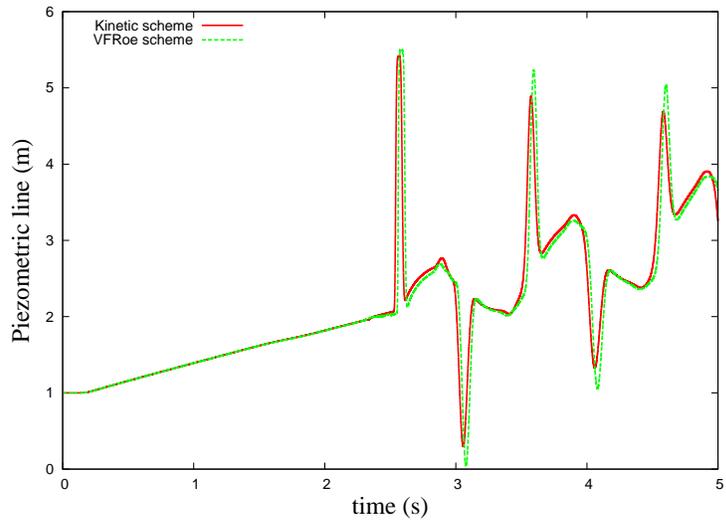}
}
\caption[Optional caption for list of figures]{Comparison between the kinetic scheme and the upwinded VFRoe scheme at $x = 0.5\, m$.}
\label{KineticVsVFRoe}
\end{center}
\end{figure}

The result is in a good agreement and is represented on \fig\ref{KineticVsVFRoe}.
\paragraph{Upwinding of the friction.\newline} 
It is well-known that cell-centered approximation of source terms leads to, 
generally, wrong results. We consider the kinetic scheme with the upwinded friction and the cell-centered one 
(i.e. we use $\WW=(Z,S,\cos\theta)$ instead of \eqref{WW} and we add the cell-centered friction $K_i^n u_i^n  \abs{u_i^n }$ 
to the right hand side of Equation \eqref{FVScheme}). We compare the schemes in a symmetrical flow.

\begin{figure}[H]
\begin{center}
\subfigure[ Cell-centered friction]
{
\includegraphics[scale = 0.85]{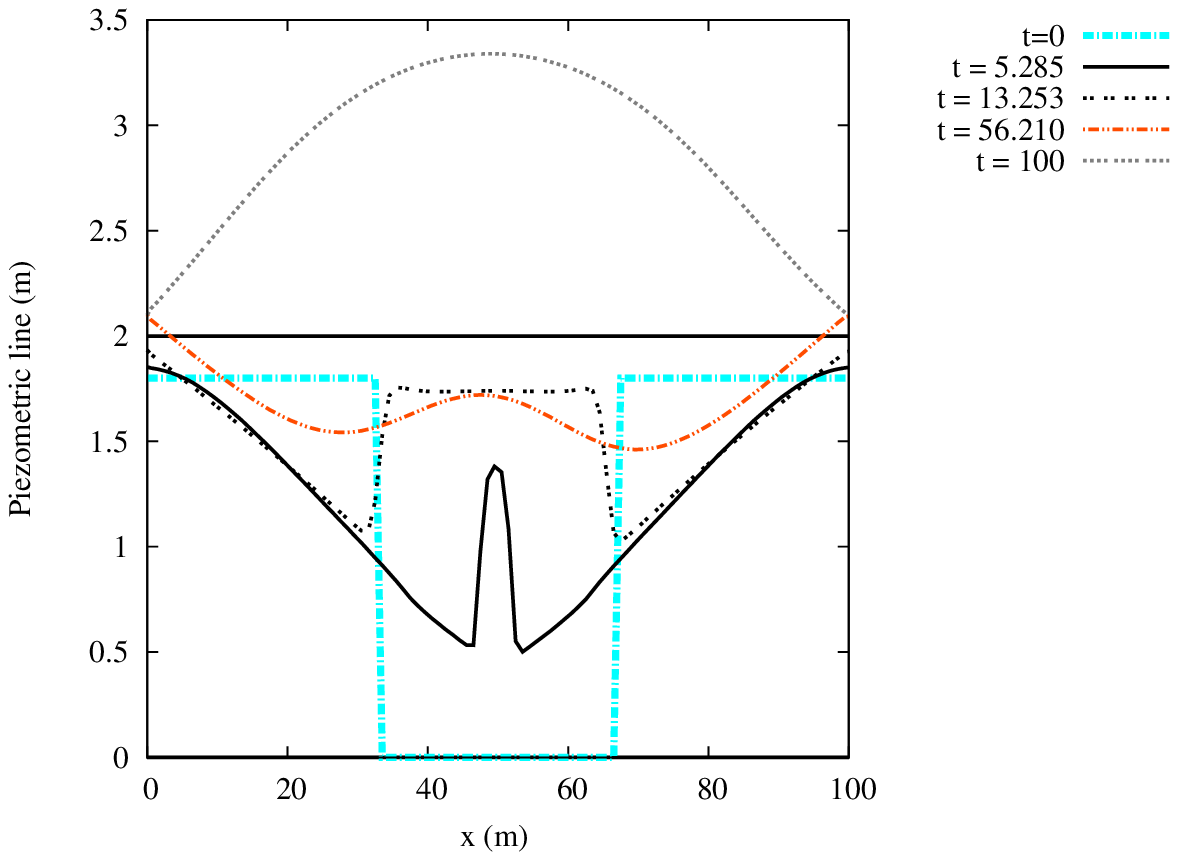}
}
\\
\subfigure[ Upwinded friction]
{
  \includegraphics[scale = 0.85]{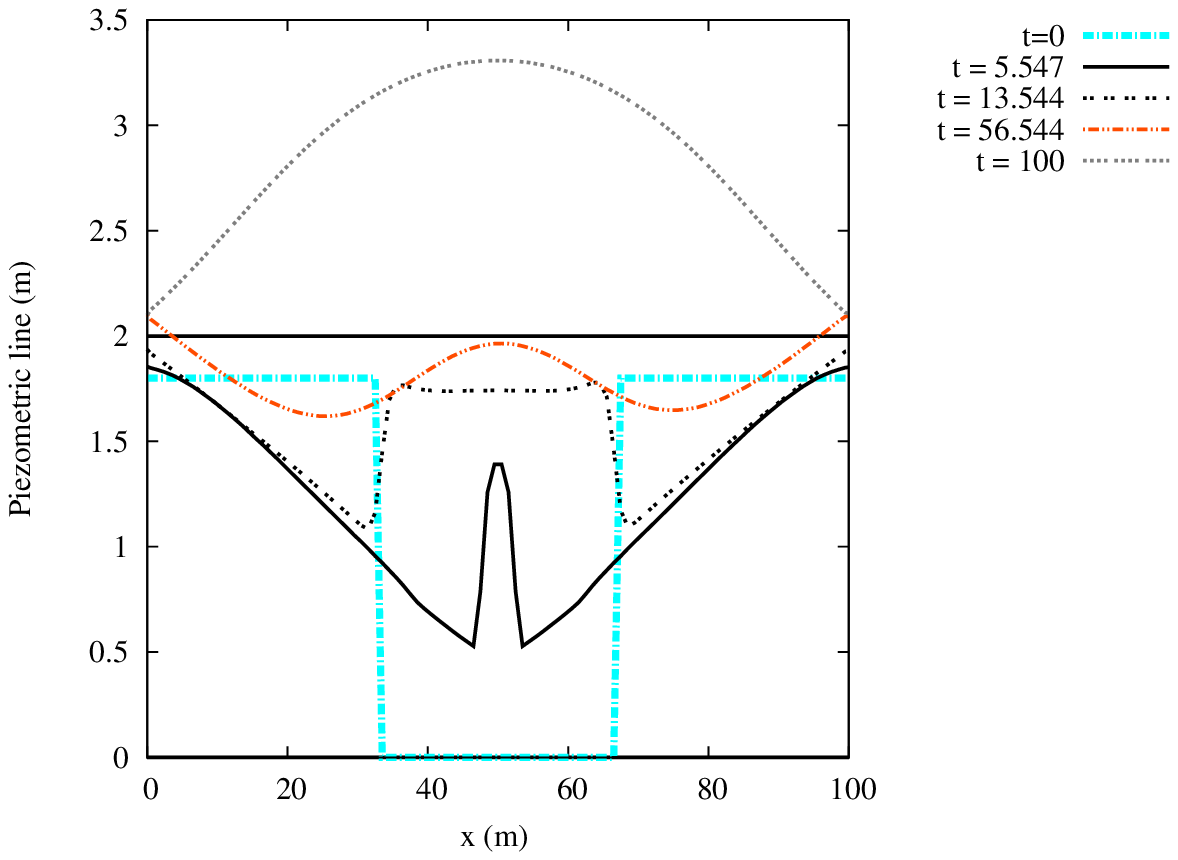}
}
\caption[Optional caption for list of figures]{Comparison of the cell-centered friction  and upwinded friction for $K_s = 100$.}
\label{FrictionKs100}
\end{center}
\end{figure}

\begin{figure}[H]
\begin{center}
\subfigure[ Cell-centered friction]
{
\includegraphics[scale = 0.85]{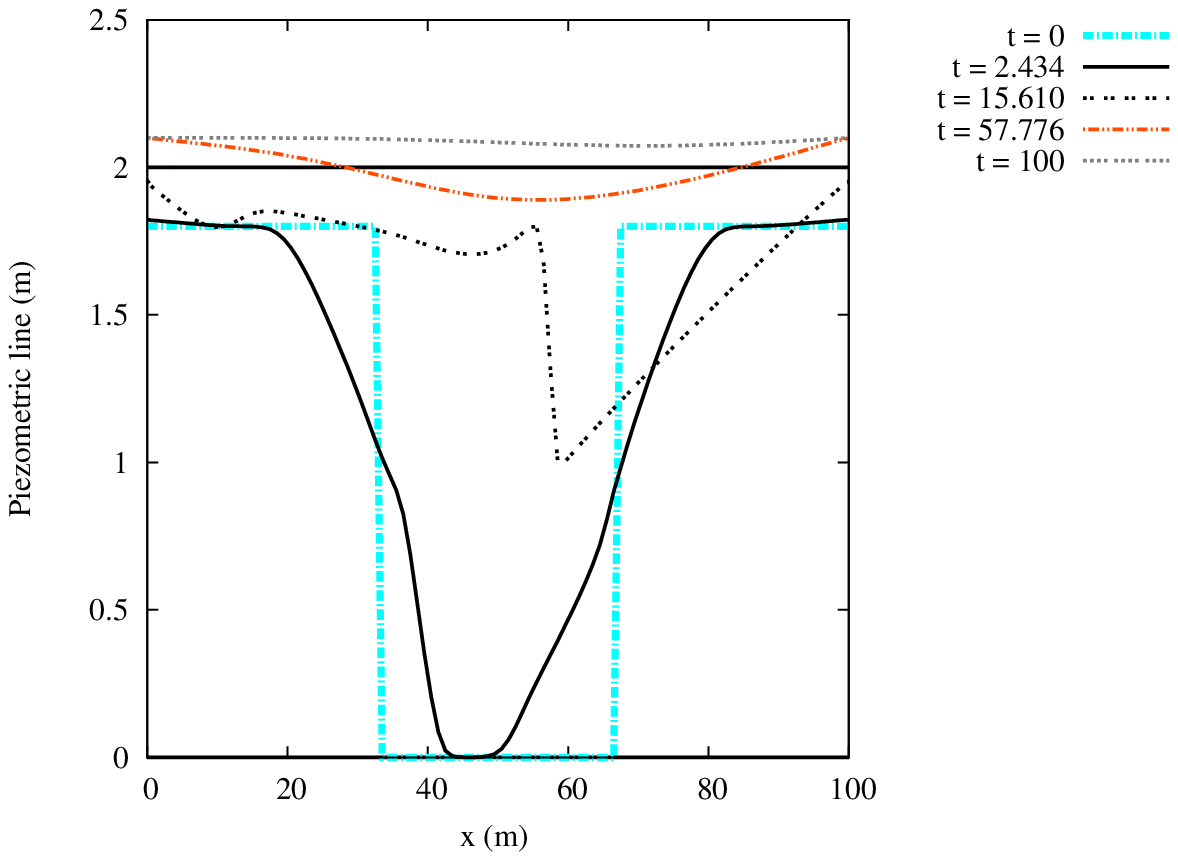}
}
\\
\subfigure[ Upwinded friction]
{
  \includegraphics[scale = 0.85]{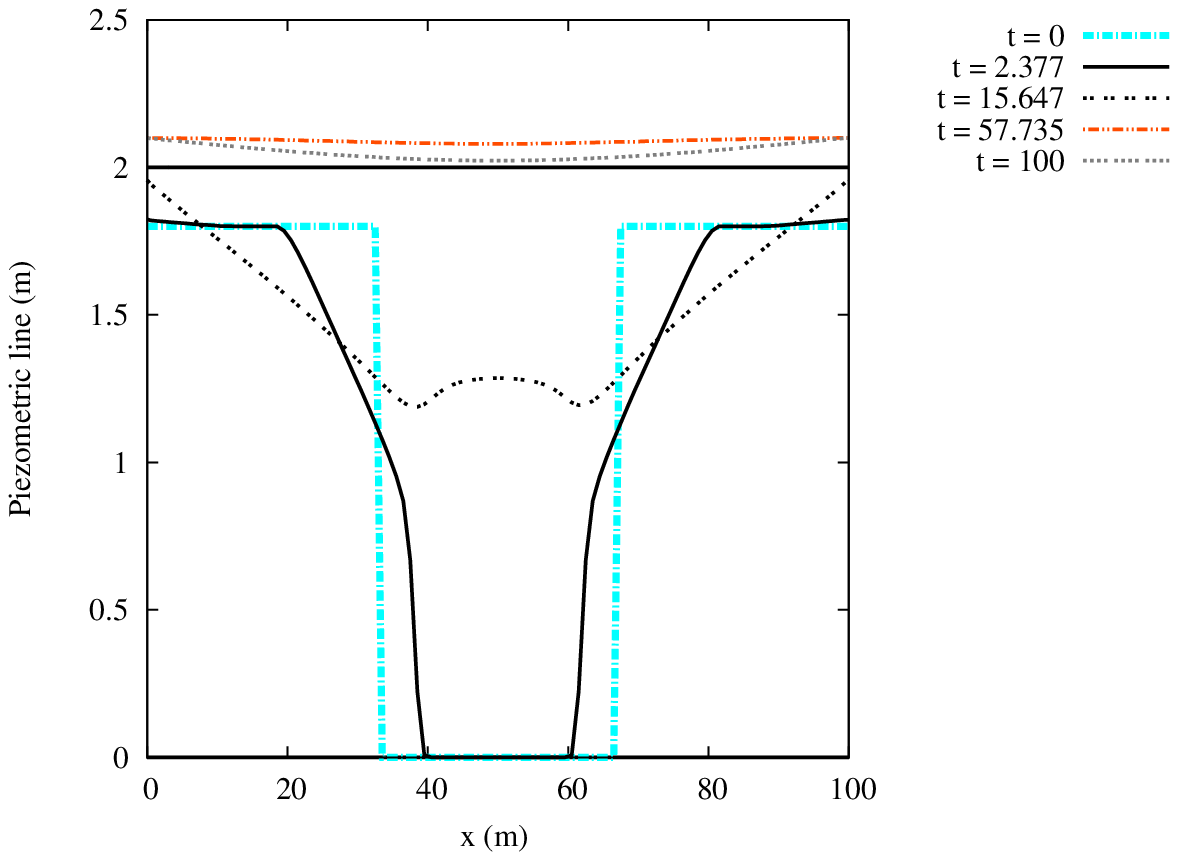}
}
\caption[Optional caption for list of figures]{Comparison of the cell-centered friction  and upwinded friction for $K_s = 10$.}
\label{FrictionKs10}
\end{center}
\end{figure}

The numerical experiment is performed on a  $100$ $m$ long closed pipe with constant section of diameter $2\;m$. 
The simulation starts from a ``double dam break'', 
as displayed on \fig\ref{FrictionKs100} and \fig\ref{FrictionKs10} at time $t=0$.
The upstream and downstream condition are identical: the
piezometric head increases linearly from $1$ to $2.1$ meters. 
We choose the same parameters as in the previous experiment.

The  results in \fig\ref{FrictionKs100}-\ref{FrictionKs10} show that the scheme with the cell-centered friction, 
contrary to the upwinded one, does not preserve the symmetry of
the flow. In particular, for $K_s=100$ (low friction) and  at time $t=56.210$ (see
\fig\ref{FrictionKs100} on top) we observe a 
small disymmetry,  which evolves drastically at time $t=1.095$ for $K_s=10$ (high friction) 
(see \fig\ref{FrictionKs10} on top). Despite the unavailability of experimental data, the kinetic
scheme with the upwinded friction term, from a physical point of view, gives the expected result,
namely, a symmetrical flow. 
\section{Conclusion}
We have presented a global manner to upwind conservative and non conservative source terms. 
To this end, we have used the definition of the non-conservative product of \cite{DLM95} 
which allows to recover the classical upwinding of conservative terms. 
Using the notion of  \emph{dynamic slope}, we have also upwinded the friction term 
given by the Manning-Strickler law (which is neither conservative nor
non-conservative)  in a FVK framework. 
The combination of all these quantities into a single one  is
an elegant and easy way to construct a kinetic scheme with reflections by introducing the potential
bareer.
Although kinetic schemes naturally deal with drying and flooding areas, the friction term is
manually set to $0$ 
when such cells appear.
\small{
\bibliographystyle{plain}

}

\end{document}